\documentclass[12pt]{article}
\usepackage{latexsym}

\def\bs{\begin{subequations}}
\def\es{\end{subequations}}

\catcode`\@=11

\newtoks\@stequation
\def\subequations{\refstepcounter{equation}
  \edef\@savedequation{\the\c@equation}%
  \@stequation=\expandafter{\theequation}
  \edef\@savedtheequation{\the\@stequation}
  \edef\oldtheequation{\theequation}%
  \setcounter{equation}{0}%
  \def\theequation{\oldtheequation\alph{equation}}}

\def\endsubequations{\setcounter{equation}{\@savedequation}%
  \@stequation=\expandafter{\@savedtheequation}%
  \edef\theequation{\the\@stequation}\global\@ignoretrue}

\makeatletter
        \renewcommand{\theequation}{\thesection.\arabic{equation}}%
        \@addtoreset{equation}{section}%
\makeatother

\renewcommand{\thefootnote}{\fnsymbol{footnote}}

\begin{document}

\begin{titlepage}
December 4, 2008

\begin{center}        \hfill   \\
            \hfill     \\
                                \hfill   \\

\vskip .25in

{\large \bf Calculus with a Quaternionic Variable \\}

\vskip 0.3in

Charles Schwartz\footnote{E-mail: schwartz@physics.berkeley.edu}

\vskip 0.15in

{\em Department of Physics,
     University of California\\
     Berkeley, California 94720}
        
\end{center}

\vskip .3in

\vfill

\begin{abstract}

Most of theoretical physics is based on the mathematics of functions 
of a real or a complex variable; yet we frequently are drawn to try 
extending our reach to  include quaternions. The 
non-commutativity of the quaternion algebra poses obstacles for the 
usual manipulations of calculus; but we show in this paper how  
many of those obstacles can be overcome.  The surprising result is that
the first order 
term in the expansion of $F(x+\delta)$ is a compact formula 
involving both $F^{\prime}(x)$  
and $[F(x) - F(x^{*})]/(x-x^{*})$. This advance in the differential 
calculus for quaternionic variables also leads us to some progress in 
studying integration.

\end{abstract}

\vfill

\end{titlepage}

\renewcommand{\thefootnote}{\arabic{footnote}}
\setcounter{footnote}{0}
\renewcommand{\thepage}{\arabic{page}}
\setcounter{page}{1}

\section{Introduction}
We are very familiar with functions of a real or complex variable $x$ 
which we can expand, in the mode of differential calculus, as
\begin{equation}
F(x+\delta) = F(x) + F^{\prime}(x) \delta + \frac {1}{2}F^{\prime 
\prime}(x) \delta^{2} + \ldots.\label{1.1}
\end{equation}

But what if we consider a quaternionic variable
\begin{equation}
x = x_{0} + i x_{1} + j x_{2} + k x_{3} ,\label{1.2}
\end{equation}
involving four real variables, $x_{\mu},\;\mu=0,1,2,3$, along with  
 those quaternions $i, j,k$ which do not commute with one another.
\begin{equation}
i^{2}=j^{2}=k^{2} = -1, \;\;\; ij=-ji=k, \;\;etc.\label{1.3}
\end{equation}
The small quantity $\delta$ will also involve all those quaternions. How 
then can we expect anything as neatly packaged as Equation (\ref{1.1})?

This is a long-standing challenge to mathematicians. The earliest attempt  
to extend the usual concept of the derivative, $\frac{dF}{dx}$, 
with a quaternionic $dx$, failed dramatically. The subsequent approach focused 
on the four real variables,
\begin{equation}
dF(x) = \sum_{\mu}\;\frac{\partial F(x)}{\partial x_{\mu}} 
dx_{\mu}.\label{1.4}
\end{equation}
That approach, often called quaternionic analyticity, springs from the 
work in the 1930's by R. Fueter \cite{1} and his school, with more accessible 
articles reviewing that subject available in references \cite{2}, \cite{3}.
Some more recent attempts to advance that work may be found in 
references \cite{4}, \cite{5}, \cite{6}.  In Appendix C I provide a 
rough summary of the Fueter approach.

The first new result presented in this paper is an alternative approach 
to the differential calculus - something between relying on the whole 
quaternionic variable, $\frac{dF}{dx}$,  and resorting to the four-component
real variables, as in Eq. (\ref{1.4}).
This starts, in Section 2, with the separation of the quaternionic displacement 
$\delta$ into two parts, one ``parallel'' and the other 
``perpendicular'' to the quaternionic variable $x$ as it may be 
envisioned in that four-dimensional space.

The subsequent Sections show how this leads to a surprisingly compact 
and general formula for the quaternionic version of the expansion 
(\ref{1.1}):
\begin{equation}
F(x+\delta) = F(x) + F^{\prime}(x)\;\delta_{\parallel} + 
(F(x) - F(x^{*}))/(x-x^{*})\;\delta_{\perp} + O(\delta^{2}).\label{1.5}
\end{equation}

In the second part of this paper we look at integration; and find 
that the new form of the quaternionic differential leads to new 
results in this other realm of calculus. 

\section{Local Coordinates}
The standard approach to quaternionic variables starts with a global set 
of imaginary coordinates;
\begin{equation}
x = x_{0} + i x_{1} + j x_{2} + k x_{3}, \;\;\;\;\; \delta = 
\delta_{0} + i \delta_{1} + j \delta_{2} + k \delta_{3}.  \label{2.1}
\end{equation}
We now want to write $x$ in a different way:
\begin{equation}
x = x_{0} + u_{x}\;r; \;\;\;\;\; r = 
\sqrt{x_{1}^{2}+x_{2}^{2}+x_{3}^{2}}, \;\;\;\;\; u_{x}^{2} = -1
\label{2.2}
\end{equation}
where $u_{x}$ is a unit imaginary that varies in its $i,j,k$ 
composition as $x$ moves from one point to another in that 
4-dimensional space. This is analogous to the use of polar 
coordinates in 2-dimensional Euclidean space.

Now we want to decompose the quaternionic quantity $\delta$ in  a 
particular way that refers to this local coordinate system.
\begin{equation}
\delta=\delta_{\parallel} + \delta _{\perp}, \;\;\;\;\; 
\delta_{\parallel} = \frac{1}{2}(\delta - u_{x}\;\delta\;u_{x}), 
\;\;\;\;\; \delta_{\perp} = \frac{1}{2}(\delta + 
u_{x}\;\delta\;u_{x}),\label{2.3}
\end{equation}
which leads to the algebraic relations
\begin{equation}
 \delta_{\parallel}\; u_{x} = u_{x}\;\delta_{\parallel}, \;\;\;\;\; 
\delta_{\perp}\;u_{x} = - u_{x}\;\delta_{\perp}. \label{2.4}
\end{equation}
The essence of this approach is expressed in the nomenclatures  
``parallel'' and ``perpendicular'' for these two components of 
$\delta$ as they relate to the local  quaternion $x$.
The most useful way to write these relations is
\begin{equation}
 \delta_{\parallel}\; x = x\; \delta_{\parallel}, \;\;\;\;\;
\delta_{\perp}\; x = x^{*}\; \delta_{\perp}, \label{2.5}
\end{equation}
where $^{*}$ is the complex conjugation operator that changes the sign 
of all imaginaries. \footnote{The efficacy of this technique was 
discovered as the result of a more long-winded calculation, which may 
be seen in reference \cite{7}.}

Now we shall give three examples of how to expand $F(x+\delta)$ 
with this simple machinery.

\section{The Function $F(x) = x^{n}$}

We calculate directly,
\begin{equation}
(x+\delta)^{n} = x^{n} + 
\sum_{m=0}^{n-1}\;x^{n-m-1}\;\delta\;x^{m} + O(\delta^{2}).\label{3.1}
\end{equation}
Putting in the separation $\delta= \delta_{\parallel} + 
\delta_{\perp}$ 
and using the properties of Eq.(\ref{2.5}), the sum becomes
\begin{equation}
\sum = \sum_{m=0}^{n-1}( x^{n-1}\;\delta_{\parallel} + 
x^{n-m-1}\;x^{*\;m}\;\delta_{\perp}) = \\
n\;x^{n-1}\;\delta_{\parallel} + 
(x^{n}-x^{*\;n})(x-x^{*})^{-1}\;\delta_{\perp},\label{3.2}
\end{equation}
where we evaluated a finite geometric series.

\section{The Exponential Function}

For a general quaternionic variable x, we can define the exponential 
function in the usual way:
\begin{equation}
e^{x} = \lim_{N\rightarrow \infty} (1+\frac{x}{N})^{N} \label{4.1}
\end{equation}
and this leads us to the expansion, 
\begin{equation}
e^{(x+\delta)} = e^{x}[1 + \int_{0}^{1}ds\; e^{-sx}\; \delta\;e^{sx} 
+ O(\delta^{2})],\label{4.2}
\end{equation}
which is correct for the situation where $x$ and $\delta$ do not commute. For a 
derivation of this formula, see Appendix A.

Putting in a real parameter $p$, and following the course set above, 
we get the expansion, 
\begin{eqnarray}
e^{p(x+\delta)} - e^{px} = \int_{0}^{1}ds\; p \;
e^{(1-s)px}(\delta_{\parallel}+ \delta_{\perp}) e^{spx} = \\ 
\int_{0}^{1} ds\; p\; e^{px}\;\delta_{\parallel} + \int_{0}^{1}ds\; 
p\; 
e^{(1-s)px}\;e^{spx*}\;\delta_{\perp} = \\
p\;e^{px}\;\delta_{\parallel} + (e^{px} - 
e^{px*})(x-x^{*})^{-1}\;\delta_{\perp}\label{4.3}
\end{eqnarray}
to first order in $\delta$.

\section{General Analytic Function F(x)}

For a general analytic function $F(x)$ of a quaternionic variable $x$, we 
start by assuming a representation as a Laplace transform:
\begin{equation}
F(x) = \int dp\; f(p)\; e^{px}\label{5.1}
\end{equation}
where $p$ is a real variable.  We then use the result of the previous 
Section to obtain
 \begin{equation}
F(x+\delta) - F(x) = F^{\prime}(x)\;\delta_{\parallel} +
  (F(x) - F(x^{*}))\;(x - x^{*})^{-1} \; \delta_{\perp} + 
  O(\delta^{2}), 
  \label{5.2} 
\end {equation}
where $F^{\prime}(x)$ is the derivative of the function $F(x)$ 
calculated as if $x$ were a real variable.

This is our general result.  The particular result of Section 3, for 
$F(x)=x^{n}$, also fits  this general formula; and thus it also 
works for any power series $F(x) = \sum_{n}\;c_{n}\;x^{n}$.

The authors of reference \cite{4} have taken an approach somewhat  
similar to what is done here. They 
introduced a local unit imaginary (which they call \emph{iota}) that is 
the same as what we have defined as $u_{x}$. However, they limit 
their differentiations to displacements that are restricted to the  
two-dimensional space of what we call 
$\delta_{\parallel}$ without allowing any of $\delta_{\perp}$.
In that way they merely 
reproduce what is known about ordinary complex variables.

\section{Further Exercises}
Let us  define the first-order differential operator ${\cal{D}}$, from 
Eq.(\ref{5.2}), as
\begin{equation}
F(x+\delta) = F(x) + {\cal{D}}\;F(x) + O(\delta^{2})\label{6.1}
\end{equation}
with
\begin{equation}
{\cal{D}}\;F(x) = F^{\prime}(x)\;\delta_{\parallel} + 
(F(x) - F(x^{*}))\;(x - x^{*})^{-1} \; \delta_{\perp}.\label{6.2}
\end{equation}
Several interesting exercises are now suggested.

Calculate ${\cal{D}}\;(F(x)G(x))$ and verify the applicability of Leibnitz'  
rule.

Calculate $\;{\cal{D}}\; \frac{1}{G(x)}\;$ and also $\;{\cal{D}}\;F(G(x))\;$. 

\section{Alternative Arrangements}

Still another way to represent our result for the first-order 
differential is in terms of some \emph{partial} derivatives, defined as 
follows.
\begin{eqnarray}
   {\cal{D}}\;F(x) = \frac{\partial F(x)}{\partial 
 x_{\parallel}}\;dx_{\parallel} + \frac{\partial F(x)}{\partial 
 x_{\perp}}\;dx_{\perp}, \;\;\;\;\;\;\;\;\;\; \label{d9} \\ 
  \frac{\partial F(x)}{\partial x_{\parallel}} \equiv F^{\prime}(x), 
  \;\;\;\;\;\frac{\partial F(x)}{\partial x_{\perp}} \equiv 
  (F(x)-F(x^{*}))(x-x^{*})^{-1} .\label{7.1}
 \end{eqnarray}
 
 Suppose we restrict the functions $F(x)$ to be real: that is, the 
 coefficients $c_{n}$ in $F=\sum_{n}\;c_{n}\;x^{n}$ or the amplitudes 
 $f(p)$ in the Laplace transform are real numbers.  Then it is noted 
 that the terms in Eq.(\ref{5.2})  can be written with the 
 displacement quaternions, $\delta_{\parallel}$ and $\delta_{\perp}$, 
 written either to the right or to the left of their accompanying 
 factors. This is obvious in the case of $F^{\prime}(x)$, since 
 $\delta_{\parallel}$ commutes with $x$. For the second term, we know 
 that $\delta_{\perp}$ does not commute with $x$; it takes the 
 complex conjugate.  But we note that the whole 
 expression$(F-F^{*})/(x-x^{*})$ is real; therefore this 
 rearrangement is possible.
 
 The same rearrangement can be done with Eq.(\ref{d9}).
 
 These considerations lead us to note that the second term in the 
 equation for ${\cal{D}}F(x)$ can be written in terms of commutators as
 \begin{equation}
(F(x) - F(x^{*}))(x-x^{*})^{-1}\;\delta_{\perp} = [C,F(x)],\label{7.2}
\end{equation}
where $C$ is defined by
\begin{equation}
	[C,x] = \delta_{\perp}, \;\;\;\;\;C = 
	\frac{1}{x^{*}-x}\;\delta_{\perp}.\label{7.3}
\end{equation}
What is somewhat surprising about this alternative arrangement, 
Eq.(\ref{7.2}), is that the expression on the left hand side is 
manifestly non-local, involving things evaluated at the point $x$ 
and also at the remote point $x^{*}$; but the right hand side appears 
to be local, involving only $x$. This confusion is removed when one 
recognizes that $C$ is a non-local operator, involving 
$\delta_{\perp}$, which changes $x$ to $x^{*}$.

\section{Second Order Terms}

Let's return to the exponential function (\ref{4.1}) and proceed with 
 the expansion,
\begin{equation}
e^{(x+\delta)} = e^{x}[1+\int_{0}^{1}ds\; e^{-sx}\delta e^{sx} +
\int_{0}^{1}dt \int_{0} ^{1-t}ds \; e^{-(s+t)x}\delta e^{tx}\delta 
e^{sx} + O(\delta^{3})].\label{8.1}
\end{equation}

The best approach is to combine the exponential function and the 
Laplace transform from the beginning. Writing $F(x+\delta) = F(x) + 
F^{(1)} + F^{(2)} + \ldots$, we now look at
\begin{equation}
F^{(2)} =\int dp f(p)\; e^{px} p^{2}\;
\int_{0}^{1}dt \int_{0} ^{1-t}ds \; e^{-(s+t)px}\delta e^{tpx}\delta 
e^{spx}.\label{8.2}
\end{equation}
Again, we decompose $\delta$ and after a bit more work arrive at the 
result for the second order term,
\begin{eqnarray}
F^{(2)} = \frac{1}{2} F^{\prime \prime}(x)\;\delta_{\parallel}^{2} + 
(F(x)-F(x^{*}))\;(x-x^{*})^{-2}\;(\delta_{\perp}\delta_{\parallel} - \delta 
\delta_{\perp}) + \nonumber \\ 
F^{\prime}(x)\;(x-x^{*})^{-1}\;\delta\delta_{\perp} +
F^{\prime}(x^{*})\;(x^{*}-x)^{-1}\;\delta_{\perp}\delta_{\parallel}.\label{8.3}
\end{eqnarray}

It is also true that, with the first order term given as $F^{(1)}(x) 
= {\cal{D}}\;F(x)$, the second order result can be written as
\begin{equation}
F^{(2)}(x) = \frac{1}{2}\; {\cal{D}\;\cal{D}}\; F(x).\label{8.4}
\end{equation}
To verify this one needs the preliminary formulas,
\begin{equation}
  {\cal{D}}\; x = \delta, \;\;\;\;\; {\cal{D}}\; x^{*} = \delta^{*}, \;\;\;\;\;
  {\cal{D}}\; \delta = 0, \;\;\;\;\;
   {\cal{D}}\;u_{x} = \frac{1}{r}\; \delta_{\perp},\label{8.5}
\end{equation}
along with $\delta_{\parallel}^{*}\;\delta_{\perp} = 
\delta_{\perp}\;\delta_{\parallel}$ and $\delta_{\perp}^{*} = 
-\delta_{\perp}$. See further in Appendix D.

\section{Discussion on Differentials}

It is surprising how simple and how general the new results obtained here are. 
It is also noteworthy that our differential operators are no longer local: 
they involve $F(x^{*})$ along with $F(x)$.

One may ask what restrictions there are on the functions $F(x)$ 
considered above. At first, one would say that they should be real 
analytic 
functions; having terms like $xax$ where $a$ is a general quaternion 
would certainly cause trouble.\footnote{This use of the term ``real 
analytic'' differs from that found in reference \cite{3}.} One can extend this condition slightly 
by allowing $F(x)$ (but not the function $G(x)$ in Section 6) 
to be a real function with arbitrary quaternions 
multiplying from the left. That is, the power series form $F = 
\sum_{n}\;c_{n}\;x^{n}$ could have arbitrary numbers $c_{n}$.

This bias to the left-hand side can be reversed if we change the 
original steps (\ref{4.2}), setting $s \rightarrow 1-s$,  and 
(\ref{5.1}), putting $f(p)$ on the right-hand side.

The Taylor series we have discussed above are expansions about the 
origin $x=0$. In the usual complex analysis such power series may be 
about any fixed point $x=x_{f}$; but such a quaternion constant put in 
the middle of our expressions would appear to cause trouble. That 
trouble could be avoided by limiting $x_{f}$ to be real; but there 
is a better way. If we define a new quaternionic variable $y = 
x-x_{f}$ then we may proceed as done above only using the appropriate 
unit imaginary $u_{y}$, instead of the original $u_{x}$ in order to 
separate the displacement $\delta$ into ``parallel'' and 
``perpendicular'' components.

One may also ask if this general method may be applied to some other 
kind of non-commuting algebra beyond the quaternions.  I believe that 
something very similar can be done starting with a Clifford algebra.
Other examples are given in Appendix B and in reference \cite{8}.

\section{Introduction to Integration }

When the conventional approach to analyticity  of 
quaternionic functions failed in the differential calculus, the 
main push was then in the realm of integral calculus.

The key result of the Fueter school was  a third order differential equation that 
could define a ``regular'' function of a quaternionic variable, just 
as the Cauchy-Riemann equation was a first order constraint on 
functions of a complex variable $z = x+iy$. That approach is described  
roughly in Appendix C. Their result is a focus on integrals over 
a 3-dimensional surface in the 4-dimensional space.

With the constructon of the quaternionic differential operator 
${\cal{D}}$ we can do something quite different about integration, as 
is shown in the following two Sections.

\section{The Line Integral }

In ordinary calculus of functions of the real variable t, we know 
what is meant by an integral, such as $\int f(t)\;dt$. But when we 
first 
consider quaternionic (or other non-commuting) variables it is 
unclear even how to write such an expression. We shall pursue that 
path in Section 12.

Alternatively, we can start with the 
defining relation between the integral and the differential:
\begin{equation}
 \int_{a}^{b}\; df(t) = \int_{a}^{b} \frac{df(t)}{dt}\;dt=f(b)-f(a) 
 ;\label{11.1}
\end{equation}
and this is what we shall generalize for our non-commuting quaternionic 
variable $x$ as, 
\begin{equation}
\int_{a}^{b}\; {\cal{D}} F(x) = F (x_{b}) - F(x_{a}).\label{11.2}
\end{equation}

We  define this integral as an additive operation along a path in 
that four-dimensional space of the real variables $x_{\mu}$, 
\begin{equation}
x = x_{path}(s), \;\;\;\;\; x_{path}(0) = x_{a}, \;\;\;\;\; 
x_{path}(1) = x_{b}\label{11.3}
\end{equation}
where $s$ is a real continuous parameter.

Next, we  subdivide 
that path, whatever it may be, into a large number of infinitesimal 
increments. 
\begin{equation}
\int_{a}^{b} = \sum_{n=1}^{n=N} \;\int^{(n)}, \;\;\;\;\;\;\;\;   
\int^{(n) }=\int_{x_{n-1}}^{x_{n}}, \;\;\;\;\; n=1, 
\ldots,N\label{11.4}
\end{equation}
where $x_{0} = x_{a}$ and $x_{N} = x_{b}$.

In any segment of this path we choose the line of integration, with 
the integrand ${\cal{D}} F(x)$,  to be 
the sum of two infinitesimal parts: 
\begin{equation}
x_{n} - x_{n-1} = \delta = \delta_{\parallel} + 
\delta_{\perp}.\label{11.5}
\end{equation}
The first part is  ``parallel'' to 
the direction of $x$ at that point, giving the contribution
\begin{equation}
\int_{\parallel}\;{\cal{D}} F(x) = F^{\prime}(x) 
\delta_{\parallel}.\label{11.6}
\end{equation}
Then the second part is ``perpendicular'', giving the contribution
\begin{equation}
\int_{\perp}\; {\cal{D}} F(x) = [F(x)-F(x^{*})](x-x^{*})^{-1}\; 
\delta_{\perp}.\label{11.7}
\end{equation}

The sum of these two parts is thus nothing other than
\begin{equation}
F(x_{n}) - F(x_{n-1})\label{11.8}
\end{equation}
to first order in the interval $\delta$.  The entire sum  then 
  results in Eq. (\ref{11.2}).

Another general proof can proceed as follows. If we start with the 
coordinate along the path $x(s) = x_{path}(s)$, then we can simply write,
\begin{equation}
{\cal{D}} x(s) = ds \frac{dx(s)}{ds}\label{11.9}
\end{equation}
since there is no commutativity problem in this representation. It is 
also true that we can express any function composed of powers of $x$ as 
\begin{equation}
F(x(s)) = A(s) + B(s)\;x(s)\label{11.10}
\end{equation}
where $A$ and $B$ are real functions, the only quaternions being in 
the single factor $x(s)$.  We then see that the integral becomes 
quite ordinary:
\begin{equation}
\int_{a}^{b}\; {\cal{D}} F(x(s)) = \int_{0}^{1}\; ds 
\frac{dF(x(s))}{ds} = F(x(s))|_{0}^{1} = F(x_{b}) - 
F(x_{a}).\label{11.11}
\end{equation}

Since this differential operator ${\cal{D}}$ 
obeys the Leibnitz rule we get the identity, usually 
called ``integration by parts'',
\begin{equation}
\int_{a}^{b} F(x)\;{\cal{D}} G(x) = F(x_{b})G(x_{b}) - 
F(x_{a})G(x_{a}) - \int_{a}^{b}\;({\cal{D}}F(x))\;G(x). \label{11.12}
\end{equation}

Loosly speaking, integration is the inverse of differentiation.  What 
we see in Eqs. (\ref{11.1}) and (\ref{11.2}) is one statement of that 
relationship. But there is also the other form, which is stated 
for real variables as
\begin{equation}
\frac{d}{dt}\; \int^{t}\; f(t')\;dt' = f(t).\label{11.13}
\end{equation}
For our quaternionic variables we start by looking at
\begin{equation}
{\cal{D}}_{x}\; \int ^{x}\; {\cal{D}}_{x'}\;F(x') \label{11.14}
\end{equation}
and then apply the first differential operator to the coordinate $x$
in two parts: first the $\delta_{\parallel}$ part and then the 
$\delta_{\perp}$ part. The result is just the integrand evaluated at 
the point $x$:
\begin{equation}
 = {\cal{D}}_{x}\;F(x); \label{11.15}
\end{equation}
and this is just what we should expect from the right hand side of 
Eq. (\ref{11.2}), with $x_{b}$ replaced by $x$.

\section {The Other Line Integral}

If we look at the common real integral and try to guess how 
to generalize it to the non-commutative quaternions, we might start with, 
\begin{equation}
\int f(t)\;dt  \; \stackrel{?}{\longrightarrow} \;  \frac{1}{2}\;\int\; (dx\; F(x) + 
F(x)\;dx) \;;\label{12.1}
\end{equation}
but why should $dx$ only appear on the outside; why not also in the 
middle of the function $F(x)$?

Let's try a most symmetrical arrangement with the function $F(x) = 
x^{n}$:
\begin{equation}
\int t^{n}\; dt \;\stackrel{?}{\longrightarrow} \; 
\frac{1}{n+1}\;\int\;(dx\; x^{n} + 
x\;dx\;x^{n-1} + x^{2}\;dx\;x^{n-2} + \ldots + x^{n}\;dx).\label{12.2}
\end{equation}
But we can recognize that the long expression in parentheses on the 
right hand side of this is nothing other than ${\cal{D}} x^{n+1}$:
\begin{equation}
{\cal{D}} F(x) \equiv  F(x+dx) -F(x), \;\; to\; first\; order\; in \; 
dx .\label{12.3}
\end{equation}
So we would then write,
\begin{equation}
\int t^{n}\; dt\; \longrightarrow\; \frac{1}{n+1}\int {\cal{D}} 
x^{n+1} = \frac{x^{n+1}}{n+1},\label{12.4}
\end{equation}
using our defining Eq. (\ref{11.2}). Now, this looks quite familiar.

We can extend this to any power series and thus offer the following 
rule.  For any analytic function of a real variable $f(t)$, for which 
we know the integral,
\begin{equation}
\int\;f(t)\;dt = h(t),\label{12.5}
\end{equation}
we can make the correspondence to quaternionic integration as follows:
\begin{equation}
\int\;f(t)\;dt \;\longrightarrow \;\int\;{\cal{D}}h(x) = 
h(x).\label{12.6}
\end{equation}
While this may look trivial for real and complex variables, it is 
something new for non-commuting variables. This arises because we have 
carefully defined and   
studied  the operator ${\cal{D}}$.

\section{Discussion on Integration}

Following what was stated earlier, we do require the functions 
$F(x)$ to be real analytic functions along the path of integration. 
Terms such as $xax$ would be allowed only for real constants $a$. 

Our first new result Eq. (\ref{11.2}) implies that the result of the 
integration depends only on the end points and is independent of the 
path. This is true if we also require that the function $F(x)$ be 
single valued. Then, we have the result that the integral over any 
closed path, ending up at the same point where it started, is zero. 
This is a significant new result, carrying the world of contour 
integration over from the complex domain to the quaternionic.

Our second new result, Eq's. (\ref{12.5}) and (\ref{12.6}), opens up 
considerable  
possibilities for integration of quaternionic functions.

\vskip 0.5cm
\begin{center} \textbf{ACKNOWLEDGMENT} \end{center}

 I am grateful to J. Wolf for some helpful conversation.

\vskip 0.5cm
\setcounter{equation}{0}
\def\theequation{A.\arabic{equation}}
\boldmath
\noindent{\bf Appendix A - Expansion of the Exponential}
\unboldmath
\vskip 0.5cm

Here we give a derivation of the formula (\ref{4.2}) for any 
non-commuting quantities $x$ and $\delta$.
\begin{eqnarray}
e^{(x+\delta)} = \lim_{N\rightarrow \infty} [1+\frac{x}{N} + 
\frac{\delta}{N}]^{N} = \;\;\;\;\;\;\;\;\;\; \\ 
\lim_{N\rightarrow \infty}\{[1+\frac{x}{N}]^{N} +
\sum_{m=0}^{N-1}\;[1+\frac{x}{N}]^{N-m-1} \;\frac{\delta}{N}\; [1+\frac{x}{N}]^{m} + 
O(\delta^{2})\}.
\end{eqnarray}
In taking the limit $N \rightarrow \infty$, we convert the sum over 
$m$ to an integral over $s = \frac{m}{N}$ and this yields
\begin{equation}
e^{(x+\delta)} = e^{x} + \int_{0}^{1}\; ds\; e^{(1-s)x}\;\delta \; 
e^{sx} + O(\delta^{2}).
\end{equation}

\vskip 0.5cm
\setcounter{equation}{0}
\def\theequation{B.\arabic{equation}}
\boldmath
\noindent{\bf Appendix B - SU(2,C)}
\unboldmath
\vskip 0.5cm

Here we shall extend the general method used above for a quaternionic 
variable to something built on a Lie Algebra - 
specifically SU(2).

Here is the Lie algebra:
\begin{equation}
[J_{1},J_{2}]=J_{3}, \;\;\; [J_{2},J_{3}]=J_{1},\;\;\; 
[J_{3},J_{1}]=J_{2},\label{B1}
\end{equation}
where the three  $J$'s are understood to be matrices over 
the complex numbers.  In particular we shall use the relations
\begin{equation}
e^{\theta J_{3}}\;J_{1}\;e^{-\theta J_{3}} = J_{1}\;cos\theta + 
J_{2}\;sin\theta, \;\;\;\;\; 
e^{\theta J_{3}}\;J_{2}\;e^{-\theta J_{3}} = J_{2}\;cos\theta - 
J_{1}\;sin\theta,\label{B2}
\end{equation}
which follow from (\ref{B1}).

The new variable $x$ is to be constructed with four real parameters as
\begin{equation}
x = x_{0}\;I + x_{1}\;J_{1} + x_{2}\;J_{2}+ x_{3}\;J_{3} \label{B3}
\end{equation}
and we want to expand $F(x+\delta) = F(x) + F^{(1)} + O(\delta^{2})$,
where $\delta$ is a small quantity in that same space of matrices as $x$.
Our first step is to define a local coordinate system at the given 
point $x$. By a suitable linear transformation (rotation) of the Lie 
algebra we make the coordinate $x$ appear as
\begin{equation}
x = x_{0}\;I + r J_{3}\label{B4}
\end{equation}
where we recognize that $r^{2} = x_{1}^{2}+x_{2}^{2}+x_{3}^{2}$.

We can now separate the displacement $\delta = \delta_{\parallel} + 
\delta_{\perp}$ as follows.
\begin{equation}
\delta_{\parallel} = \delta_{0}\;I + \delta_{3}\;J_{3}, \;\;\;\;\;
\delta_{\perp} = \delta_{1}\;J_{1} + \delta_{2}\;J_{2}.\label{B5}
\end{equation}

Now we are ready to study the first order term in the expansion, 
again using the representation of $F(x)$ in terms of the exponential 
function.
\begin{equation}
F^{(1)} = \int dp\;f(p)\;p\;e^{px}\;\int_{0}^{1}ds\; e^{-spx}\;\delta 
\;e^{spx}.\label{B7}
\end{equation}
Since $\delta_{\parallel}$ commutes with $x$, the first part of this 
is simply $F^{\prime}(x)\;\delta_{\parallel}$.  For the part with 
$\delta_{\perp}$ we use the formulas (\ref{B2}), where $\theta$ is 
replaced by $-spr$.  The integrals over $s$ are trivial and we 
merely write $sin(pr)$ and $cos(pr)$ in terms of $e^{\pm ipr}$ to get 
our final result.
\begin{eqnarray}
F(x+\delta) = F(x) + F^{\prime}(x)\delta_{\parallel} + 
\{F(x+ir)-F(x-ir)\}\;\frac{1}{2ir}\;\delta_{\perp} + \nonumber \\
\{F(x+ir) + F(x-ir) - 2F(x)\} \;\frac{1}{2r}\; [J_{3},\delta_{\perp}] +\;\;
 O(\delta^{2}).\;\;\;\;\;\label{B8}
\end{eqnarray}

It should be noted that the $\delta$-related factors in Eq. 
(\ref{B8}) can be written in the following way:
\begin{eqnarray}
[J_{3},\delta_{\perp}] = \frac{1}{r}\;[x,\delta] , \label{B9}\\ 
\delta_{\perp} = - \frac{1}{r^{2}}\;[x,[x,\delta]],\label{B10} \\
\delta_{\parallel} = \delta - \delta_{\perp}\label{B11}.
\end{eqnarray}
This means that we do not have to carry out the ''rotation'' that 
gave us Eq. (\ref{B4}) explicitely; the talk about choosing a local 
coordinate system is merely rhetorical. 

I expect that this method can be extended to other Lie algebras, with 
the quantity $\delta_{\perp}$ subdivided into distinct portions according to 
the roots of the particular algebra. The system of Eqs. 
(\ref{B10}), (\ref{B11}) would be adapted to make those separations, 
using the known values of the roots; and those root values would also 
appear in the final generalization of Eq. (\ref{B8}). 

Extension of this method to general matrix variables, over 
the complex numbers, is given in reference \cite{8}.

\vskip 0.5cm
\setcounter{equation}{0}
\def\theequation{C.\arabic{equation}}
\boldmath
\noindent{\bf Appendix C - Fueter's Differential Equation}
\unboldmath
\vskip 0.5cm

The  literature on Fueter's analysis of quaternionic functions 
points to a third order differential equation as his key result, 
extending the familiar Cauchy-Riemann (first order) equation for functions of a complex 
variable,
\begin{equation}
(\frac{\partial}{\partial x} + i\frac{\partial}{\partial y})f(z=x+iy) 
= 0.\label{C1}
\end{equation}

Here I wish to present a rather simple derivation of that result, 
starting with the exponential function of our quaternionic variable.
\begin{eqnarray}
x =x_{0}+ ix_{1}+jx_{2}+kx_{3} = x_{0} + \textbf{u}\cdot \textbf{x}, 
\;\;\;\;\; \textbf{u} = (i,j,k),\;\;\;\;\; r^{2} = 
\textbf{x}\cdot\textbf{x}, \label{C2} \\ 
e^{px} = e^{px_{0}}(cos pr + \textbf{u}\cdot\textbf{x}\;\frac{sin 
pr}{r}).\;\;\;\;\;\;\;\;\;\;\label{C3}
\end{eqnarray}

Here is Fueter's first order differential operator.
\begin{equation}
\Box = \frac{\partial}{\partial x_{0}} + i\frac{\partial}{\partial 
x_{1}} + j\frac{\partial}{\partial x_{2}} + k\frac{\partial}{\partial 
x_{3}} = \frac{\partial}{\partial x_{0}} + \textbf{u}\cdot \nabla 
,\label{C4}
\end{equation}
and we calculate its action on the exponential function and find the 
result
\begin{equation}
\Box\;e^{px} = -2\;e^{px_{0}}\;\frac{sin pr}{r}.\label{C5}
\end{equation}
This is not zero (as in the complex case, Eq.(\ref{C1})) but it is rather 
simple (and real). Moreover, we recognize this function as a solution of the 
four-dimensional Laplace equation.
\begin{equation}
\triangle_{4} = \Box\;\Box^{*} = \frac{\partial^{2}}{\partial 
x_{0}^{2}} + \triangle_{3}.\label{C6}
\end{equation}
So here are two forms of Fueter's third order differential equation.
\begin{equation}
\triangle_{4}\;\Box\;e^{px} = 0, \;\;\;\;\; 
\triangle_{4}\;\Box\;x^{n} = 0,\label{C7}
\end{equation}
where the second result comes from expanding the first result in a 
power series in the parameter p.  Thus any superposition (with real 
coefficients) of the powers or the exponential will satisfy this 
condition; and this is the basis for what they define as holomorphic 
 functions of a quaternionic variable. They also exclude functions 
 with terms such as $xax$ with arbitrary quaternion constants $a$ - 
 just as we have done in the present paper.

 From those differential equations (\ref{C7}) (Cauchy-Riemann-Fueter), some integral theorems follow. In the case of 
 complex variables (\ref{C1}), we get the familiar result that any integral around a 
 closed path in the complex plane will be zero (with suitable 
 analytic and single-valued behavior of the function f(z)). In the 
 quaternionic case, the relevant integral  is over a closed 
 3-dimensional surface 
 in the 4-dimensional space of the $x_{\mu}$.
 
  Here is a surprise!  Look at the result of our differential 
 operator ${\cal{D}}$ acting on the exponential function $e^{px}$, 
 Eq. (\ref{4.3}). The coefficient of $\delta_{\perp}$ is the same 
 function that we see on the right hand side of Eq. (\ref{C5}). So, 
 here is a new identification:
 \begin{equation}
\frac{\partial F(x)}{\partial x_{\perp}} = -\frac{1}{2}\;\Box F(x).
\end{equation}

\vskip 0.5cm
\setcounter{equation}{0}
\def\theequation{D.\arabic{equation}}
\boldmath
\noindent{\bf Appendix D - Caculating} ${\cal{D}}\;u_{x}$
\unboldmath
\vskip 0.5cm

Here is a derivation of the last item in Eq. (\ref{8.5}), which is 
made easy if we take a geometric perspective as we write the 
coordinate in four dimensions as $x = x_{0}+ru_{x}$.  Start by 
writing 
\begin{equation}
{\cal{D}}\;u_{x} = \alpha\;\delta_{\parallel} + \beta \; 
\delta_{\perp},
\end{equation}
where $\alpha$ and $\beta$ are to be determined.  First, consider a 
displacement that has only $\delta_{\parallel}$: this should not 
change $u_{x}$ at all, since such a displacement can only change 
$x_{0}$ and $r$.  Thus, we see that $\alpha = 0$.

Next, consider a displacement that has only $\delta_{\perp}$: this 
should not change $x_{0}$ or $r$. So we write,
\begin{equation}
{\cal{D}}\;u_{x} = {\cal{D}}_{\perp}\; u_{x} = {\cal{D}}_{\perp}\; 
(x-x_{0})/r = \frac{1}{r}\;{\cal{D}}_{\perp}\;x = 
\frac{1}{r}\;\delta_{\perp}.
\end{equation}

Following the result for $F^{(2)}$  in Section 8, one may wonder 
whether the entire Taylor series might be written as,
\begin{equation}
F(x+\delta)  = 
\sum_{k=0}^{\infty}\;\frac{1}{k!}\; {\cal{D}}^{k}\;F(x) 
= e^{{\cal{D}}}\;F(x).
\end{equation}
This may be readily verified for the 
functions $F(x)=x^{n}$, starting with the equation,
\begin{equation}
e^{{\cal{D}}}\;x\;e^{-{\cal{D}}} = x + \delta. \label{D.4}
\end{equation}
For the exponential function, define 
\begin{equation}
Q(p) = e^{{\cal{D}}}\;e^{px};
\end{equation}
then calculate $dQ/dp$ and use Eq. (\ref{D.4}).


\begin{thebibliography}{99}
	
\bibitem{1} R. Fueter, {\sl Comment. Math. Helv.}\/ {\bf 7}, 307 
(1935) and {\bf 8}, 371 (1936).

\bibitem{2} C. A. Deavours, {\sl Amer. Math. Monthly}\/ {\bf 80}, 995 
(1973).

\bibitem{3} A. Sudbery, {\sl Math. Proc. Camb. Phil. Soc.}\/ {\bf 85}, 
199 (1979).

\bibitem{4}
S. De Leo and P. P. Rotelli, {\sl Appl. Math. Lett.}\/ {\bf 16}, 
1077 (2003); arXiv:funct-an/9703002 

\bibitem{5}
G. Gentili and C. Stoppato, arXiv:0802.3861 [math.CV]

\bibitem{6}
D. Alayon-Solarz, arXiv:0803.3480v2 [math.CV] 

\bibitem{7}
C. Schwartz, arXiv:0803.3782 [math.FA] 

\bibitem{8}
C. Schwartz, arXiv:0804.2869 [math.FA] 


\end{thebibliography}
\end{document}